\documentclass[12pt]{amsart}
\usepackage{amsmath}
\usepackage{amsfonts,amssymb,amsthm}
\newcommand{\thm}[2]{\begin{#1} #2 \end{#1}}

\newcommand{\btheta}{\boldsymbol{\theta}}

\newtheorem{theorem}{Theorem}[section]

\newcommand{\integers}{{\bf Z}}

\newcommand{\reals}{{\bf R}}

\begin{document}

%-------------- Author entries --------------------

\title{Symmetrized Chebyshev polynomials}
%Article title
%\shorttitle{DRAFT} % Shortened version for
                                             % headline title

\author{Igor Rivin}

%\address{Mathematics Department, California Institute of Technology,
%Pasadena, CA 91125}

\address{Department of Mathematics, Temple University, Philadelphia}

\curraddr{Mathematics Department, Princeton University}

\email{rivin@math.temple.edu}

\thanks{These results first appeared in the author's 1998 preprint
``Growth in free groups (and other stories),'' but seem to be of
independent interest. The positivity result was preprint
math.CA/0301210, but there appears to be no reason to separate it
from the limiting distribution result, and many reasons to keep
them together.}

\date{today}



\begin{abstract}
We define a class of multivariate Laurent polynomials closely
related to Chebyshev polynomials, and prove the simple but
somewhat surprising (in view of the fact that the signs of the
coefficients of the Chebyshev polynomials themselves alternate)
result that their coefficients are non-negative. We further show
that a Central Limit Theorem holds for our polynomials.
\end{abstract}

\maketitle

\section*{Introduction}

Let $T_n(x)$ and $U_m(x)$ be the Chebyshev polynomials of the
first and second kinds, respectively. We define the following
Laurent polynomials (which, for lack of a better name, we call the
\emph{symmetrized Chebyshev polynomials}:

\begin{eqnarray}
R_n(c; x_1, \dots, x_n) &=& T_n\left(\frac{c\sum_{i=1}^k \left[x_i
+ \frac{1}{x_i}\right]}{2 k}\right),\\
S_n(c; x_1, \dots, x_n) &=& U_n\left(\frac{c\sum_{i=1}^k \left[x_i
+ \frac{1}{x_i}\right]}{2 k}\right).
\end{eqnarray}
The function $R_n$ arises in the enumeration of conjugacy classes
in the free group on $k$ generators, more specifically, there is
the following result:
\begin{theorem}[\cite{rwalks}]
  \label{homoenum} The number of
cyclically reduced words of length $k$ in $F_r$ homologous to $e_1
[a_1] + \cdots + e_r [a_r]$ is equal to the coefficient of
$x_1^{e_1} \cdots x_r^{e_r}$ in
\begin{equation}
\label{genfn} 2\left(\sqrt{2r-1}\right)^k R_k({\frac{r}
{\sqrt{2r-1}}}; x_1, \dots, x_r)  + (r-1)[1 + (-1)^k]
\end{equation}
\end{theorem}

\medskip\noindent
{\em Remark.} The rescaled Chebyshev polynomial $T_k(a x)/a^k$ is
called the $k$-th Dickson polynomial $T_k(x, a)$ (see
\cite{schur}).

The coefficients of generating functions of combinatorial objects
are non-negative, so one is led to wonder for which values of $c$
are the coefficients of $R_n(c; x_1, \dots, x_k)$ and $S_n(c; x_1,
\dots, x_k)$ nonnegative. In this note we give an essentially
complete answer (see Theorems \ref{nonneg} and \ref{mnonneg}). We
also write down an explicit formula (Eq. (\ref{fullform})) for the
coefficients of $R_n$ and $S_n.$ Furthermore, we use the
positivity and tools of probability theory to analyze the
distribution of the coefficients of $R_n$ to prove a central limit
theorem -- Theorem \ref{centlim}

\section{Some facts about Chebyshev polynomials}
\label{chebint}

The literature on Chebyshev polynomials is enormous; \cite{rivlin}
is a good to start. Here, we shall supply the barest essentials in
an effort to keep this paper self-contained.

There are a number of ways to define Chebyshev polynomials (almost
as many as there are of spelling their inventor's name). A
standard definition of the {\em Chebyshev polynomial of the first
kind} $T_n(x)$ is:

\begin{equation}
\label{def1} T_n(x) = \cos n \arccos x.
\end{equation}
In particular, $T_0(x) = 1,$ $T_1(x) = x.$ Using the identity
\begin{equation}
\label{cosid} \cos(x+y) + \cos(x-y) = 2 \cos x \cos y
\end{equation}
we immediately find the three-term recurrence for Chebyshev
polynomials:
\begin{equation}
\label{threerec} T_{n+1}(x) = 2 x T_n(x) - T_{n-1}(x).
\end{equation}
The definition of Eq. (\ref{def1}) can be used to give a ``closed
form'' used in Section \ref{limitu}:
\begin{equation}
\label{sqrtdef} T_n(x) = {\frac12}\left[\left(x -
\sqrt{x^2-1}\right)^n + \left(x + \sqrt{x^2-1}\right)^n\right].
\end{equation}
Indeed, let $x = \cos \theta.$ Then
$$\left(x -
\sqrt{x^2-1}\right)^n = \exp(-i n \theta),$$ while
$$\left(x +
\sqrt{x^2-1}\right)^n = \exp(i n \theta),$$ so that
$${\frac12}\left(x - \sqrt{x^2-1}\right)^n + \left(x +
\sqrt{x^2-1}\right)^n = \Re \exp(i n \theta) = \cos n \theta.$$

We also define Chebyshev polynomials of the second kind $U_n(x)$,
which can again be defined in a number of ways, one of which is:
\begin{equation}
\label{derivdef} U_n(x) = {\frac1{n+1}}T_{n+1}^\prime(x).
\end{equation}
A simple manipulation shows that if we set $x = \cos \theta,$ as
before, then
\begin{equation}
\label{trigdef} U_n(x) = \frac{\sin (n+1) \theta}{\sin \theta}.
\end{equation}
In some ways, Schur's notation $\mathcal{ U}_n = U_{n-1}$ is
preferable. In any case, we have $U_0(x) = 1$, $U_1(x) = 2 x,$ and
otherwise the $U_n$ satisfy the same recurrence as the $T_n$, to
wit,
\begin{equation}
\label{trec2} U_{n+1}(x) = 2 x U_n(x) - U_{n-1}(x).
\end{equation}
From the recurrences, it is clear that for $f=T, U$, $f_n(-x) =
(-1)^n f(x),$ or, in other words, every second coefficient of
$T_n(x)$ and $U_n(x)$ vanishes. The remaining coefficients
alternate in sign; here is the explicit formula for the
coefficient $c_{n-2m}^{(n)}$ of $x^{n-2m}$ of $T_n(x):$
\begin{equation}
\label{coefform} c_{n-2m}^{(n)} = (-1)^m {\frac{n}{n-m}}
{\binom{n-m}{m}} 2^{n-2m-1}, \qquad m=0, 1, \ldots,
\left[{\frac{n}2}\right].
\end{equation}
This can be proved easily using Eq. (\ref{threerec}).

\section{Analysis of the functions $R_n$ and $S_n$.}
\label{genanal}

In view of the alternation of the coefficients, the appearance of
the Chebyshev polynomials as generating functions in Section
\ref{homoesec} seems a bit surprising, since combinatorial
generating functions have non-negative coefficients. Below we
state and prove a generalization. Remarkably, Theorems
\ref{nonneg} and \ref{mnonneg} do not seem to have been previously
noted.

\thm{theorem} { \label{nonneg} Let $c > 1.$ Then all the
coefficients of $R_n(c; x)$ are non-negative. Indeed the
coefficients of $x^n, x^{n-2}, \ldots, x^{-n + 2}, x^{-n}$ are
positive, while the other coefficients are zero. The same is true
of $S_n$ in place of $R_n.$ }

\begin{proof} Let $a_n^k$ be the coefficient of $x^k$ in
$U_n((c/2)(x+1/x)).$ The recurrence gives the following recurrence
for the $a_n^k:$
\begin{equation}
\label{newrec} a_{n+1}^k = c(a_n^{k-1}+a_n^{k+1}) - a_{n-1}^k.
\end{equation}
Now we shall show that the following always holds:

\begin{description}
\item[(a)] $a_n^k \geq 0$ (inequality being strict if and only if
$n-k$ is even).

\item[(b)] $a_n^k \geq \max(a_{n-1}^{k-1}, a_{n-1}^{k+1}),$ the
inequality strict, again, if and only if $n-k$ is even.

\item[(c)] $a_n^k \geq a_{n-2}^k$ (strictness as above).
\end{description}

The proof proceeds routinely by induction; first the induction
step (we assume throughout that $n-k$ is even; all the quantities
involved are obviously $0$ otherwise):

By induction $a_{n-1}^k < \min(a_n^{k-1}, a_n^{k+1}),$ so by the
recurrence \ref{newrec} it follows that $a_{n+1}^k >
\max(a_n^{k-1}, a_n^{k+1}).$ (a) and (c) follow immediately.

For the base case, we note that $a_0^0 = 1,$ while $a_1^1 =
a_1^{-1} = c > 1,$ and so the result for $U_n$ follows. Notice
that the above proof does {\em not} work for $T_n$, since the base
case fails. Indeed, if $b_n^k$ is the coefficient of $x^k$ in
$T_n((c/2)(x+1/x))$, then $b_0^0 = 1$, while $b_1^1 = c/2$, not
necessarily bigger than one. However, we can use the result for
$U_n$, together with the observation (which follows easily from
the addition formula for $\sin$) that
\begin{equation}
\label{moretrig} T_n(x) = {\frac{U_n(x) - U_{n-2}(x)}2}.
\end{equation}
Eq. (\ref{moretrig}) implies that $b_n^k = a_n^k - a_n^{k-2} > 0$,
by (c) above.
\end{proof}

The proof above goes through almost verbatim to show:

\thm{theorem} { \label{mnonneg} Let $c > 1.$ Then all the
coefficients of $R_n$ are non-negative. The same is true of $S_n$
in place of $R_n$ }

To complete the picture, we note that:

\thm{theorem} { \label{trivthm}
$$
R_n(1; x) = {\frac12}\left(x^n + {\frac1 {x^n}}\right).
$$
}

\begin{proof} Let $x = \exp i\theta.$ Then $1/2(x+1/x) = \cos \theta,$
and $R_n(1; x)=T_n(1/2(x+1/x)) = \cos n \theta = 1/2(x^n +
1/x^n).$
\end{proof}

\thm{remark} { For $c<-1$ it is true that all the coefficients of
$R_n(c; .)$ and $S_n(c; .)$ have the same sign, but the sign is
$(-1)^n.$ For $|c|<1,$ the result is completely false. For $c$
imaginary, the result is true. I am not sure what happens for
general complex $c$. }

By the formula (\ref{coefform}), we can write
\begin{equation}
\label{uniexp} T_n\left({\frac{c}
2}\left(x+{\frac{1}{x}}\right)\right) = {\frac12}
\sum_{m=0}^{\left[{\frac{n}{2}}\right]} (-1)^m {\frac{n}{n-m}}
{\binom{n-m}{m}} c^{n-2m} \left(x+{\frac{1}{x}}\right)^{n - 2m}.
\end{equation}
Noting that
\begin{equation}
\label{binom} \left(x+{\frac{1}{x}}\right)^k = \sum_{i=0}^k
{\binom{k} { i}} x^{k - 2i}
\end{equation}
we obtain the expansion
\begin{equation}
\label{fullform} R_n(c; x) = c^n \sum_{k=-n}^n x^k
\sum_{m=0}^{\left[{\frac{n}{2}}\right]}
\left(-{\frac{1}{c^2}}\right)^m {\frac{n}{n-m}} {\binom{n-m}{m}}
{\binom{n-2m}{(n-2m-k)/2}},
\end{equation}
where it is understood that $\binom{a} { b}$ is $0$ if $b<0,$ or
$b > a$, or $b \notin \integers.$ A similar formula for $S_n$ can
be written down by using Eq. (\ref{moretrig}).

\section{Limiting distribution of coefficients}
\label{limitu}

While the formula (\ref{fullform}) is completely explicit, and a
similar (though more cumbersome) expression could be obtained for
$R_n(c; x_1, \dots, x_k),$ for many purposes it is more useful to
have a limiting distribution formula as given by Theorem
\ref{centlim} below. To set up the framework, we note that since
all the coefficients of $R_n(c; x_1, \dots, x_k)$ are non-negative
(according to Theorem \ref{mnonneg}), they can be thought of
defining a probability distribution on the integer lattice
$\integers^k,$ defined by $p(l_1, \dots,
l_k)=[x_1^{l_1}x_2^{l_2}\cdots x_k^{l_k}]R(c; x_1, \dots,
x_k)/R(c; 1, \dots, 1)$ (where the square brackets mean that we
are extracting the coefficients of the bracketed monomial). Call
the resulting probability distribution $\mathcal{ P}_n(c; {\bf
z}),$ where ${\bf z}$ now denotes a $k$-dimensional vector.

\thm{theorem} { \label{centlim} With notation as above, when
$c>1$, the probability distributions $\mathcal{ P}_n(c; {\bf
z}/\sqrt{n})$ converge to a normal distribution on $\reals^k$,
whose mean is ${\bf 0}$, and whose covariance matrix $C$ is
diagonal, with entries
$$
\sigma^2 = {\frac{c}{k}}\left[1+
\left({\frac{c+1}{c-1}}\right)^{1/2}\right].
$$
}

To prove Theorem \ref{centlim} we will use the method of
characteristic functions (Fourier transforms), and more
specifically at first the {\em Continuity Theorem} (\cite[Chapter
XV.3, Theorem 2]{feller}), \thm{theorem} { \label{ct} In order
that a sequence $\{F_n\}$ of probability distributions converges
properly to a probability distribution $F$, it is necessary and
sufficient that the sequence $\{\phi_n\}$ of their characteristic
functions converges pointwise to a limit $\phi$, and that $\phi$
is continuous in some neighborhood of the origin.

In this case $\phi$ is the characteristic function of $F$. (Hence
$\phi$ is continuous everywhere and the convergence
$\phi_n\rightarrow \phi$ is uniform on compact sets). }

The characteristic function $\phi_n$ of $\mathcal{ P}_n(c; {\bf
z})$ is simply $$R_n(c; \exp( i \theta_1), \ldots, \exp(i
\theta_k))/R(c; 1, \ldots, 1).$$  By definition of $R_n$,

\begin{eqnarray*}
R_n(c; \exp( i \theta_1), \ldots, \exp(i
\theta_k)) &= T_n\left({\frac{c}{k}}\sum_{j=1}^k \cos \theta_j\right),\\
R_n(c; 1, \ldots, 1)) &= T_n\left({\frac{c}{k}}\sum_{j=1}^k \cos
0\right) = T_n(c).
\end{eqnarray*}

We now use the form of Eq. (\ref{sqrtdef}):

$$
T_n(x) = {\frac12}\left(\left(x - \sqrt{x^2-1}\right)^n + \left(x
+ \sqrt{x^2-1}\right)^n\right),
$$
setting
$$
u = \sum_{j=1}^k \cos {\frac{\theta_j} {\sqrt{n}}},
$$
and
$$\btheta = (\theta_1, \dots, \theta_k),$$
we get

\begin{equation}
\label{phiform} \phi_n\left({\frac{\btheta}{\sqrt{n}}}\right) =
\frac{\left\{ \left({\frac{c}{k}} u + \sqrt{{\frac{c^2}{k^2}}
u^2-1}\right)^n + \left({\frac{c}{k}} u - \sqrt{{\frac{c^2}{k^2}}
u^2-1}\right)^n\right\}}{2T_n(c)}.
\end{equation}
Notice, however, that for $c>1$, the ratio of the second term in
braces to the first is exponentially small as $n\rightarrow
\infty$, since the first term grows like $(c+\sqrt{c^2-1})^n$,
while the second as $(c-\sqrt{c^2-1})^n$ (since $\cos
{\frac{\theta_j}{\sqrt{n}}} \rightarrow 1$). Since, for the same
reason, $2T_n(c) = (c+\sqrt{c^2-1})^n[1 + o(1)],$ we can write:
$$
\phi_n({\frac{\btheta}{\sqrt{n}}}) = \left[\frac{{\frac{c}{k}} u +
\sqrt{{\frac{c^2}{k^2}} u^2-1}}{c+\sqrt{c^2-1}}\right]^n + o(1).
$$
Substituting the Taylor expansions for the cosine terms (hidden in
$u$ for typesetting reasons), we get:
\begin{equation}
u = k + {\frac1{2n}} \langle\btheta, \btheta\rangle + o(1/n),
\end{equation}
so
\begin{equation}
\label{lintrm} {\frac{c}{k}} u = c + {\frac{c}{2kn}} \langle
\btheta, \btheta\rangle + o(1/n).
\end{equation}
A similar computation gives
\begin{equation}
{\frac{c^2}{k^2}} u^2 = c^2 + {\frac{c^2}{k n}} \langle \btheta,
\btheta\rangle + o(1/n).
\end{equation}
Substituting the last expansion into the square root, we see that
\begin{eqnarray*}
\sqrt{{\frac{c^2}{k^2}} u^2 - 1} &= \sqrt{c^2-1}\sqrt{1+
{\frac1n}\left[{\frac{c^2} {(c^2-1) k}} \langle\btheta,
\btheta\rangle + o({\frac1{n}})\right]}\\
&= \sqrt{c^2-1} \left[1+{\frac1{2n}}{\frac{c^2} {(c^2-1) k}}
\langle\btheta, \btheta\rangle \right] + o({\frac1{n}}).
\end{eqnarray*}
Adding Eq. (\ref{lintrm}) and collecting terms, get
\begin{eqnarray*}
\lefteqn{\frac{c u + \sqrt{c^2 u^2-k^2}}{k
\left(c+\sqrt{c^2-1}\right)} =} \\
& & 1+{\frac1{2 n k}} \left(1+{\frac1{c+ \sqrt{c^2-1}}}\right)
\left(c + {\frac{c^2}{(c^2-1)^{1/2} }}\right) \langle \btheta,
\btheta\rangle + o({\frac1{n}}).
\end{eqnarray*}
Performing some further simplifications, we see that
$$
\phi_n({\frac{\btheta}{\sqrt{n}}}) =
\exp\left(-\frac12\btheta^\perp C \btheta\right) + o(1),
$$
where $C$ is the covariance matrix described in the statement of
Theorem \ref{centlim}, and Theorem \ref{centlim} follows
immediately.

\bibliographystyle{amsplain}

\end{document}